\documentclass[12pt,oneside]{article}
\usepackage{geometry}
\usepackage{graphicx}

\usepackage{soul}

\usepackage{xcolor}
\usepackage[
    colorlinks=true,
    urlcolor=blue,
    linkcolor=blue,
    citecolor=blue,
    filecolor=blue,
]{hyperref}

\usepackage{ifxetex}
\ifxetex
  \usepackage{fontspec}
\else
  \usepackage[utf8]{inputenc}
\fi
\usepackage{bbold}
\usepackage{graphicx}
\usepackage{url}

\usepackage{amssymb}
\usepackage{amsmath}
\usepackage{amsthm}

\newtheorem{theorem}{Theorem}
\newtheorem{definition}[theorem]{Definition}

\usepackage{biblatex}
\usepackage{sidecap}
\usepackage{float}
\usepackage{listings}

\newfloat{lstfloat}{htbp}{lop}
\floatname{lstfloat}{Listing}
\DeclareMathOperator*\img{img}

\title{3d-printing Identification Spaces of the Square}
\author{
  {Mikael Vejdemo-Johansson}
}
\date{2019-08-29}

\begin{document}

\maketitle

\begin{abstract}
  We describe three identification spaces of the square, interesting choices of immersion into $\mathbb{R}^3$, and a process to construct 3d-printable models of their parametrizations.

\resizebox{25pc}{!}{
 \includegraphics[width=0.25\linewidth]{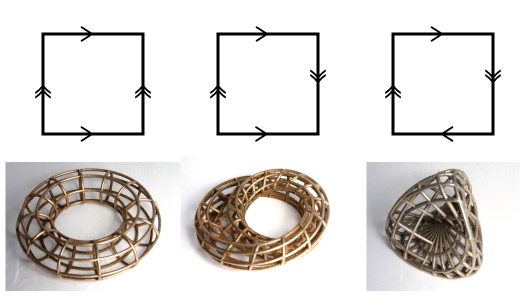}
}
\end{abstract}

The square has three \emph{identification spaces} formed by identifying opposite sides of the square.
They differ by how many of the side pairs are flipped when gluing:
\begin{enumerate}
\item[0] pair flipped: the torus
\item[1] pair flipped: the Klein bottle
\item[2] pairs flipped: the projective plane 
\end{enumerate}

These form one of the simplest examples of the \emph{fundamental polygon} construction of closed surfaces: each closed surface can be constructed by gluing sides of a polygon.

There are many ways of representing each of these as 3-dimensional objects.
For the torus, these representations can be topologically faithful -- whereas neither the Klein bottle nor the projective plane can be immersed into $\mathbb{R}^3$.
Both require self-intersections in any representation.

\begin{SCfigure}
  \centering
  \includegraphics[width=0.4\textwidth]{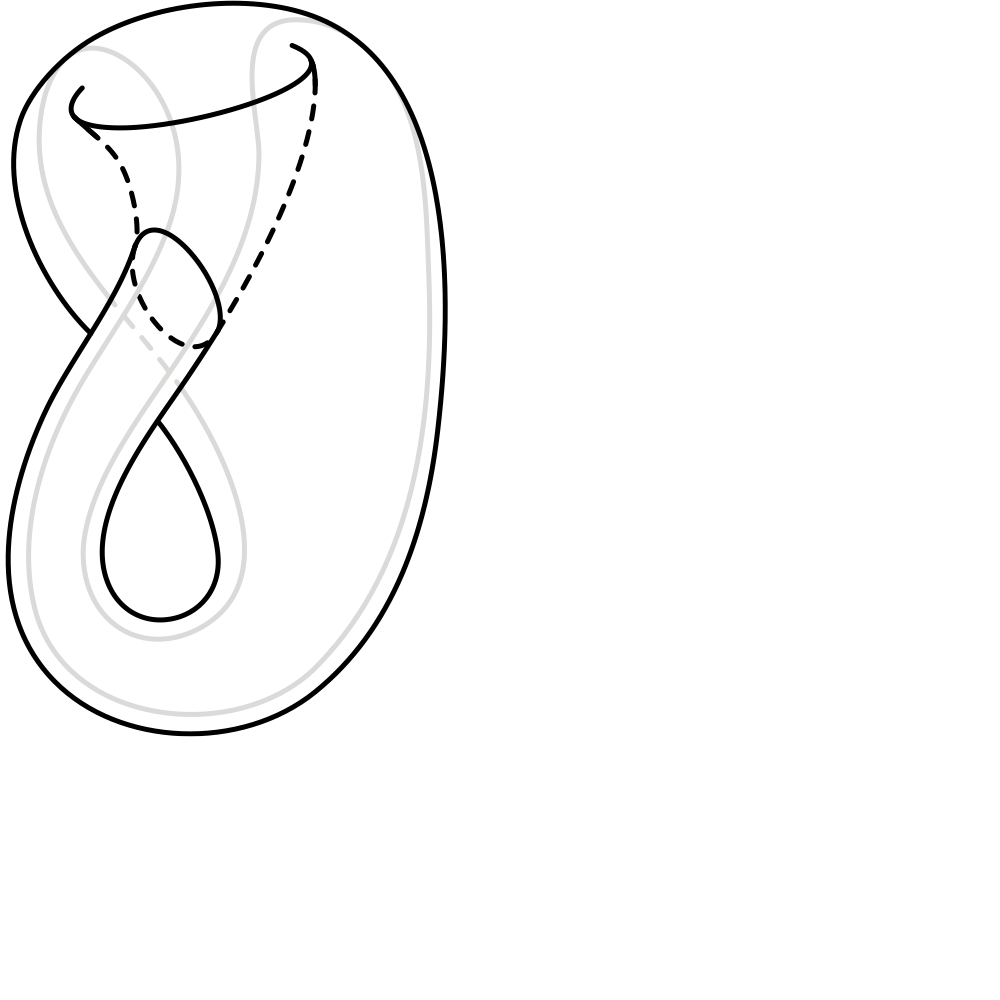}
  \caption{The most commonly used shape for the Klein bottle in $\mathbb{R}^3$. 
Picture in public domain, created by Jens Bossaert.}
  \label{fig:klein-swanneck}
\end{SCfigure}

For the Klein bottle, the most common 3d model, widely available as glass bottles or as 3d-printed objects, is one that highlights one way to view the construction: as a cylinder that self-intersects to glue the two boundary circles together.
An example can be seen in Figure~\ref{fig:klein-swanneck}.

In this paper, I will describe a 3d-printed triptych of these three identification spaces, with a different choice of Klein bottle representation than the one in Figure~\ref{fig:klein-swanneck}.
Instead I am using a projection that highlights the structure of a \emph{circle bundle}.
Circle bundles, or more generally fiber bundles, are important tools in algebraic topology that allow the use of \emph{long exact sequence} constructions when studying homotopy groups.

In Section~\emph{Fiber bundles in algebraic topology (\autoref{scrivauto:12})} I will describe the underlying algebraic topological concepts: homology, homotopy, long exact sequences and fiber bundles.
In Section~\emph{Identification spaces as bundles (\autoref{scrivauto:22})} I describe how the Klein bottle and the torus, but not the projective plane, can be constructed as circle bundles over the circle with the same construction.
In Section~\emph{3d-printing parametrizations (\autoref{scrivauto:24})} I discuss particulars of how I generated the 3d-models for printing the artwork, and include OpenSCAD code snippets for the code that generated the models.

\section{Fiber bundles in algebraic topology}
\label{scrivauto:12}

Algebraic Topology concerns itself with measuring topological properties of spaces using algebraic tools: homology groups and linear algebra, as well as homotopy groups.
Of these tools, homotopy is the more powerful, but also by far more difficult to use in calculations.

In this paper we only have space for a very condensed summary -- we recommend the excellent book by Hatcher~\cite{hatcher_algebraic_2002} for details on all the topics we touch on here.

\subsection{Homotopy}
\label{scrivauto:14}

Poincaré~\cite{poincare_1895} introduced the study of closed loops in a space as a method of studying properties of that space.
A closed loop here is a map $f:[0,1]\to X$ such that $f(0) = f(1) = x_0$ for some chosen \emph{base point} $x_0$.
We can compose loops by traversing first one and then the other, defining 

\[
(f*g)(t) = 
\begin{cases}
f(2t) & 0 \leq t \leq \frac12 \\
g(2t-1) & \frac12 \leq t \leq 1
\end{cases}
\]

We define two such loops to be \emph{homotopic} if one can be continuously deformed into the other, in other words if there is some map $H:[0,1]\times[0,1]\to X$ called a \emph{homotopy} from $f$ to $g$ such that

\begin{enumerate}
\item $H(x,0) = f(x)$
\item $H(x,1) = g(x)$
\item $H(0,t) = H(1,t) = x_0$
\end{enumerate}

We can introduce an equivalence relation on the set $\mathcal{L}(X, x_0)$ of loops with basepoint $x_0$ by setting $f\sim g$ if $f$ and $g$ are homotopic.

We can then define the \emph{fundamental group} of $X$ to be $\pi_1(X,x_0) = \mathcal{L}(X, x_0)/\sim$.
This turns into a group using $[f]*[g] = [f*g]$ as the group operation.

Homotopy theory studies fundamental groups and their higher dimensional analogues -- where instead of closed loops, the basic objects are hyperspheres mapped into $X$, with the same type of homotopy equivalence relation.

\subsection{Homology}
\label{scrivauto:16}

To present the underlying idea of homology, we will simplify the setting a little bit: introduce a discretization from an arbitrary space to focus on simplicial complexes.
Most interesting spaces can be approximated by simplicial complexes.

An \emph{$n$-simplex} is the convex hull of the unit vectors in $\mathbb{R}^{n+1}$, or equivalently the convex hull of any set of $n+1$ points in general position.
This generalizes points, intervals, triangles and tetrahedra to higher dimensions.
The convex hull of any subset of the vertices of a simplex is a \emph{subsimplex} or \emph{face} of the simplex -- this way, each $n$-simplex is built out of a collection of $2^{n+1}-1$ simplices: ${n+1\choose k+1}$ $k$-subsimplices for each dimension $0\leq k\leq n$.

A \emph{(geometric) simplicial complex} is a set of simplices, such that all intersections are subsimplices. 
In other words, a simplicial complex is formed by gluing simplices together along their faces, without allowing any intersections other than the glue surfaces.

For a simplicial complex $\Sigma$, we define its \emph{chain groups} to be abelian groups $C_k$ with one generator for each $k$-simplex.
Chains are formal linear combinations of simplices.

In homotopy, the object that measured topological properties were the closed loops: these are going to be our objects of interest in the homological setting as well.
One way to characterize a closed loop is that it has no end points.
To formalize this, we can define a \emph{boundary operator} to take a chain to its set of end points (or boundary faces in general):

\[
\partial[s_0,\dots,s_k] = 
\sum_{j=0}^k (-1)^j [s_0,\dots,s_{j-1},s_{j+1},\dots,s_k
\]

This operator takes an edge $[s,t]$ to the formal difference $t-s$, and a path $[s,t]+[t,u]$ to $(t-s) + (u-t) = u-s$.
So paths are mapped using the boundary operator to the formal differences of their end points. 
If these end points coincide, the boundary operator maps to 0.

We can take this property: $\partial z = 0$ to define what it should mean to be a \emph{cycle}.
Some of these cycles are not surprising: it is a short algebraic argument to show that $\partial^2=0$, or in other words that boundaries of things are cycles (boundaries do not in turn have boundaries).
We can extract the \emph{surprising cycles} or the \emph{essential cycles} as a quotient group: all cycles, modulo the non-surprising cycles.
From this, we define the \emph{homology} of $\Sigma$ to be 

\[
H_k\Sigma = HC_k\Sigma = \frac{\ker\partial_k}{\img\partial_k}
\]

The basis elements of $H_k\Sigma$ correspond to essential $k$-dimensional ``bubbles'' in $\Sigma$.

\subsection{Long exact sequences}
\label{scrivauto:18}

Long exact sequence are such a powerful computational tool that they found a place in the axiomatization of homology.
A sequence, here, is a sequence of abelian groups with homomorphisms between them

\[
\dots\xrightarrow{f_{i+1}} G_i
\xrightarrow{f_{i}} G_{i-1}
\xrightarrow{f_{i-1}} G_{i-2}
\xrightarrow{f_{i-2}} \dots
\]

A sequence is a \emph{chain complex} if $\img f_{i+1}\subseteq\ker f_i$ everywhere.
This is a direct abstraction of the chain complexes used to define homology above -- the requirement $\img f_{i+1}\subseteq\ker f_i$ is enough for our definition of homology.

A sequence is \emph{exact} if $\img f_{i+1} = \ker f_i$.

Many algebraic structures can be formulated in terms of exact sequences:

\begin{itemize}
\item $f$ is injective iff $0\to X \xrightarrow{f} Y$ is exact
\item $f$ is surjective iff $X \xrightarrow{f} Y \to 0$ is exact
\item $f$ is an isomorphism iff $0\to X\xrightarrow{f} Y \to 0$ is exact
\end{itemize}

A \emph{short exact sequence} is an expression of a quotient in terms of exact sequences: $0\to A\xrightarrow{f} B\xrightarrow{g} C\to 0$ induces an isomorphism $C = B / \img{f}$.

A \emph{long exact sequence} is any exact sequence longer than a short exact sequence.
These generalize the inclusion/exclusion principle algebraically.
By strategically identifying 0 modules or 0 maps within a long exact sequence, or by splitting it into interweaved short exact sequences, many computations in homological algebra and algebraic topology become easier -- or feasible at all.

Given a short exact sequence of chain complexes $0\to C_*X\to C_*Y\to C_*Z\to 0$, there is a long exact sequence in homology

\begin{multline*}
\dots
\to H_{n+1}(X)
\to H_{n+1}(Y)
\to H_{n+1}(Z) \to \\
\to H_{n}(X)
\to H_{n}(Y)
\to H_{n}(Z) \to \\
\to H_{n-1}(X)
\to H_{n-1}(Y)
\to H_{n-1}(Z)
\to\dots
\end{multline*}

These long exact sequences allow us a lot of useful calculations.

One special case here is the Mayer-Vietoris long exact sequence.
This is the long exact sequence generated by the observation that
\[
0\to 
C_*(X\cap Y) \to
C_*X \oplus C_*Y \to
C_*(X\cup Y)
\]

By splitting the sphere into a northern and a southern cap, intersecting in a neighborhood of the equator, we could use the Mayer-Vietoris sequence. Here, $X$ and $Y$ are both disks, with $H_0X = H_0Y = \mathbb{Z}$ and no other homology while the equator is a cylinder, with $H_1(X\cap Y) = H_0(X\cap Y) = \mathbb{Z}$ and no other homology.

The long exact sequence then is
\begin{multline*}
H_2X\oplus H_2Y \to H_2S^2 \to 
H_1(X\cap Y)\to H_1X \oplus H_1Y \to H_1S^2\to \\
H_0(X\cap Y)\to H_0X \oplus H_0Y \to H_0S^2 \to 0
\end{multline*}

Now $H_2X, H_2Y, H_1X, H_1Y$ all vanish, leaving us the long exact sequence
\[
0 \to H_2S^2 \to 
\mathbb{Z} \to
0 \to H_1S^2 \to
\mathbb{Z} \to 
\mathbb{Z}\oplus\mathbb{Z} \to 
H_0S^2 \to 0
\]

The subsequence $0\to H_2S^2\to \mathbb{Z}\to 0$ gives an isomorphism $H_2S^2 = \mathbb{Z}$.

The map $H_1S^2\to\mathbb{Z}$ comes from the most difficult part of the long exact sequence construction.
A cycle in the union $X\cup Y$ can be subdivided into parts that lie only in $X$ and only in $Y$, so that $z = x+y$. Then $\partial z = \partial x + \partial y = 0$ since $z$ was chosen to be a cycle. Hence $\partial x = -\partial y$, and this join boundary is in the intersection. So we map $[z]\to[\partial x]$.
Any such boundary is constructed out of paired points: $\partial x = \sum (s_i-t_i)$ -- but since the equator band is connected, each point is homologous to each other point. Thus, $[\partial x] = 0$, which means that for $H_1S^2$ the sequence reduces to $0\to 0\to H_1S^2 \to 0$, so $H_1S^2 = 0$.

This leaves the short exact sequence $0\to\mathbb{Z}\xrightarrow{x\mapsto(x,x)}\mathbb{Z}\oplus\mathbb{Z}\xrightarrow{(x,y)\mapsto x-y} H_0S^2\to 0$.
Because the sequence is exact, $H_0S^2 = \img((x,y)\mapsto x-y) = \mathbb{Z}$.

\subsection{Fiber bundles}
\label{scrivauto:20}

In general these long exact sequence arguments do not work to calculate homotopy groups. 
A \emph{fiber bundle} is a short exact sequence of spaces $0\to F\to E\xrightarrow{p} B\to 0$ with enough extra structure that long exact sequences do work.
The trick is to make the preimage subspaces $p^{-1}\subset E$ all homeomorphic to each other.
We call $F$ the \emph{fiber}, $E$ the \emph{total space} and $B$ the \emph{base space}.

One obvious source of fiber bundles are products: $0\to X\to X\times Y\xrightarrow{p} Y\to 0$ certainly has the property that $p^{-1}(y) = X\times\{y\}$ is essentially the same space for every value of $y$.
Other fiber bundles are similar to products of spaces, but can admit twists.
As a first example, two fiber bundles over the circle are the cylinder and the Klein bottle.
By using the center circle of each as the base space, the fiber over any one point is a closed interval, say $[-1,1]$.
The difference between them can only be seen globally.

More formally, following~\cite{hatcher_algebraic_2002}, we have the following important definitions and results:

\begin{definition}
  A \emph{fiber bundle structure} on a total space $E$ with fiber $F$ and base space $B$ is a projection map $p:E\to B$ such that for every $b\in B$ there is some neighborhood $b\in U$ and a homeomorphism $h: p^{-1}(U)\to U\times F$ so that $\pi_U\circ h = p$. 
\end{definition}

This last condition forces the first component of $h$ to be just $p$ itself.
The homeomorphism $h$ is called a \emph{local trivialization}, and makes each small slice of the total space look like a product space. of the local neighborhood with the fiber space.

\begin{definition}
  A map $p:E\to B$ has the \emph{homotopy lifting property} with respect to another space $X$ if given a homotopy $H:X\times [0,1]\to B$ from $f$ to $g$ in the base space and a map $\hat f:X\to E$ that lifts $f$ to $E$ -- so that $p\circ\hat f = f$ -- there is a homotopy $\hat H: X\times [0,1]\to E$ that lifts all of $H$: $p\times\mathbb{1}\circ\hat H = H$.

  A map $p:E\to B$ that has the homotopy lifting property with respect to all spaces $X$ is a \emph{fibration}.
\end{definition}

\begin{theorem}
  Suppose $p:E\to B$ has the homotopy lifting property with respect to all disks $D^k$.
  Pick basepoints $b_0\in B$ and $x_0\in p^{-1}(b_0)$.
  If $B$ is path connected, then there is a long exact sequence of homotopy groups

\begin{multline*}
\dots
\to \pi_{n+1}(F, x_0)
\to \pi_{n+1}(E, x_0)
\to \pi_{n+1}(B, b_0) \to \\
\to \pi_{n}(F, x_0)
\to \pi_{n}(E, x_0)
\to \pi_{n}(B, b_0) \to \\
\to \pi_{n-1}(F, x_0)
\to \pi_{n-1}(E, x_0)
\to \pi_{n-1}(B, b_0)
\to\dots
\end{multline*}
\end{theorem}

This ties together with a theorem generalizing work by Huebsch and Hurewicz:

\begin{theorem}
  A fiber bundle $p:E\to B$ has the homotopy lifting property with respect to all disks.
\end{theorem}

In practice this means that if we know any two of $F, E, B$ well enough we can derive much if not most of the homotopy structure for the third one.

\section{Identification spaces as bundles}
\label{scrivauto:22}

For the three spaces of interest to us -- the torus, the Klein bottle and the projective plane -- two have a clear structure of circle bundles over the circle:

\[
0\to S^1\to T^2\to S^1\to 0
\qquad
0\to S^1\to K\to S^1\to 0
\]

where the construction can be easily seen in the identification diagram: the projection $p$ collapses the square to its central line.
With the identification on the sides, this line is a circle, and preimages of small intervals -- because the top and bottom are not flipped -- are all cylinders.

For the projective plane, the situation is worse. We could still project onto a central circle, but with the flip, any preimage comes out not as a cylinder but as a Möbius strip.
The construction of the torus and the Klein bottle as circle bundles over the circle does not carry to the projective plane.
Instead, the projective plane occurs as the base space in several interesting fiber bundles.

This is what my choice of designs for the 3d-printed identification spaces is meant to illustrate.
For the torus, it is easy to see the circle bundle structure -- the usual way to depict or sculpt a torus shows the transversal circles as fibers over any chosen longitudinal circle.

For the Klein bottle, the most common immersion\footnote{Mapping such that whenever the image self-intersects, the intersecting areas have different derivatives.} into 3-dimensional space is the Swan Neck immersion.
One reason for its popularity is that it shows the construction of a Klein bottle by taking a cylinder and gluing the end points quite vividly.
For a circle bundle illustration, however, the \emph{figure 8} immersion is more appropriate: sweeping a figure 8 (seen as a circle with a self-intersection) along a Möbius strip produces the Klein bottle depicted in the visual abstract on the first page.

For the projective plane, finally, I chose to 3d-print an parametrization that does not illustrate any circle bundle structure, but instead is related to the Veronese embedding from algebraic geometry: the Roman surface.

Each of the surfaces are 3d-printed as a wireframe grid, so that the parametrizations of each can be easily seen and studied.
For the Klein bottle, some effort was added to dodge the figure 8 wires to avoid self-intersections in the fibers.
This highlights both the self-intersection structure of the surface and avoids the more confusing shape of having self-intersecting fibers: the circles in the fibers are quite clear to see.

\section{3d-printing parametrizations}
\label{scrivauto:24}

These 3d-prints were produced using OpenSCAD~\cite{openscad} for modeling, and with Shapeways for production.
For a successful 3d-print, one fundamental requirement is for the model to be \emph{watertight}

%%% MAYBE CITE HENRY SEGERMAN'S BOOK FOR THIS?
-- in other words that the mesh representation is a union of closed oriented surfaces, with no gaps.
Components may, however overlap freely.

To simplify smoothing, and construction of curves in a wireframe representation of a parametrization of a surface, I decided to express each model as a union of small watertight ``capsules'': convex hulls of a pair of spheres placed closely together.

I place these capsules along the gridlines of a square grid mapped through a parametrization $f(u,v)$ of the surface itself.
For each edge $(u,v) \to (u+1,v)$ I place capsules along $f(u+k\epsilon,v) \to f(u+(k+1)\epsilon, v)$ and for each edge $(u,v) \to (u,v+1)$ I place capsules along $f(u,v+k\epsilon) \to f(u,v+(k+1)\epsilon)$.
These sequences trace out the geodesic wires along the wireframe.

In OpenSCAD, this process comes down to two primary building blocks.
First, a function representation of the parametrization $f(u,v)$, and secondly a core loop that places all these capsules.
The three chosen parametrizations can be seen in code listings \ref{lst:torus}, \ref{lst:klein} and \ref{lst:roman}.

\textbf{Torus}
For the Torus, I used an explicit construction. Through 
\[
u\mapsto (R+r\cos u, 0, r\sin u)
\]
and an explicit rotation matrix, rotating by an angle $v$ around the axis $(0,0,1)$, the point $(u,v)\in[0,360]^2$ parametrizes a point on a torus with inner radius $r$ and outer radius $R$.
Glossing over the implementation of a rotation matrix, the resulting OpenSCAD code for the parametrization can be found in Listing~\ref{lst:torus}.

\begin{lstfloat}
  \begin{lstlisting}
function torus(i,j) = 
  let(u = i*360/lat_ribs)
  let(v = j*360/long_ribs)
  let(ix = inner_radius*cos(u)) 
  let(iy = 0)
  let(iz = inner_radius*sin(u))
  [ix+outer_radius,iy,iz]*rotate_matrix(0,0,v);
  \end{lstlisting}
\caption{Parametrization code for the Torus embedding.}
\label{lst:torus}
\end{lstfloat}

\begin{figure}
\begin{center}
  \includegraphics[width=0.49\textwidth]{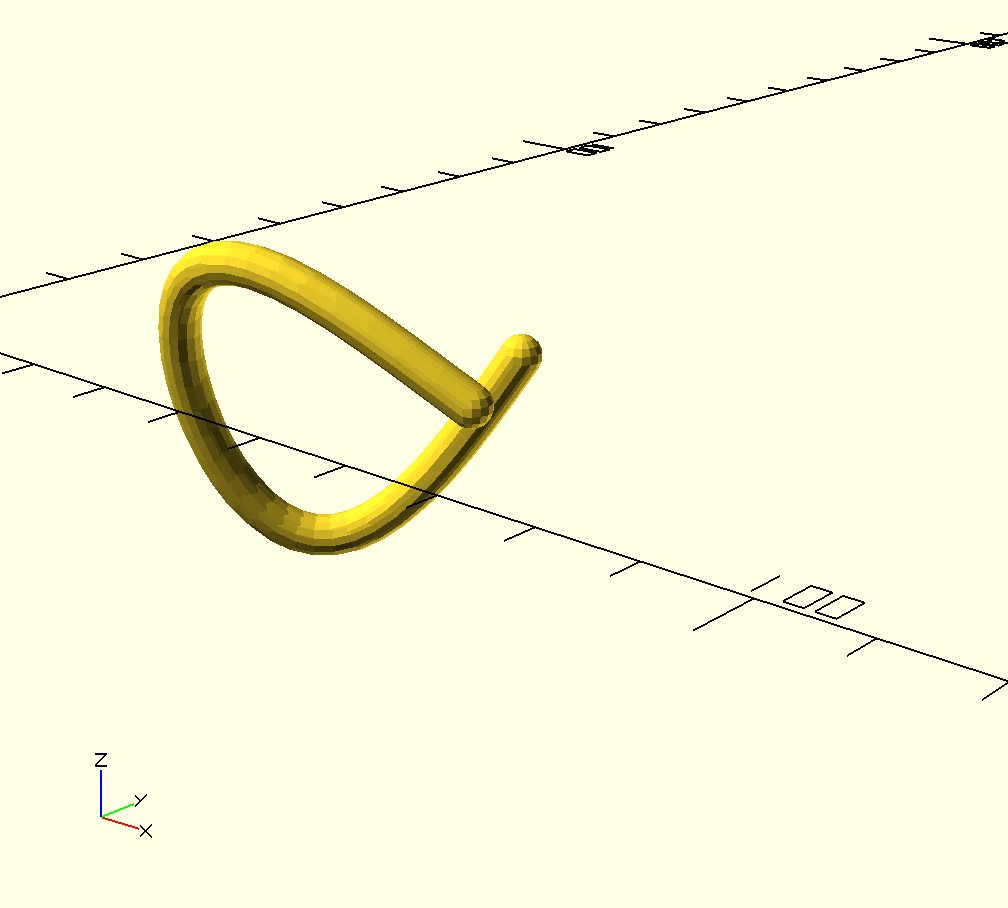}\hfill
  \includegraphics[width=0.49\textwidth]{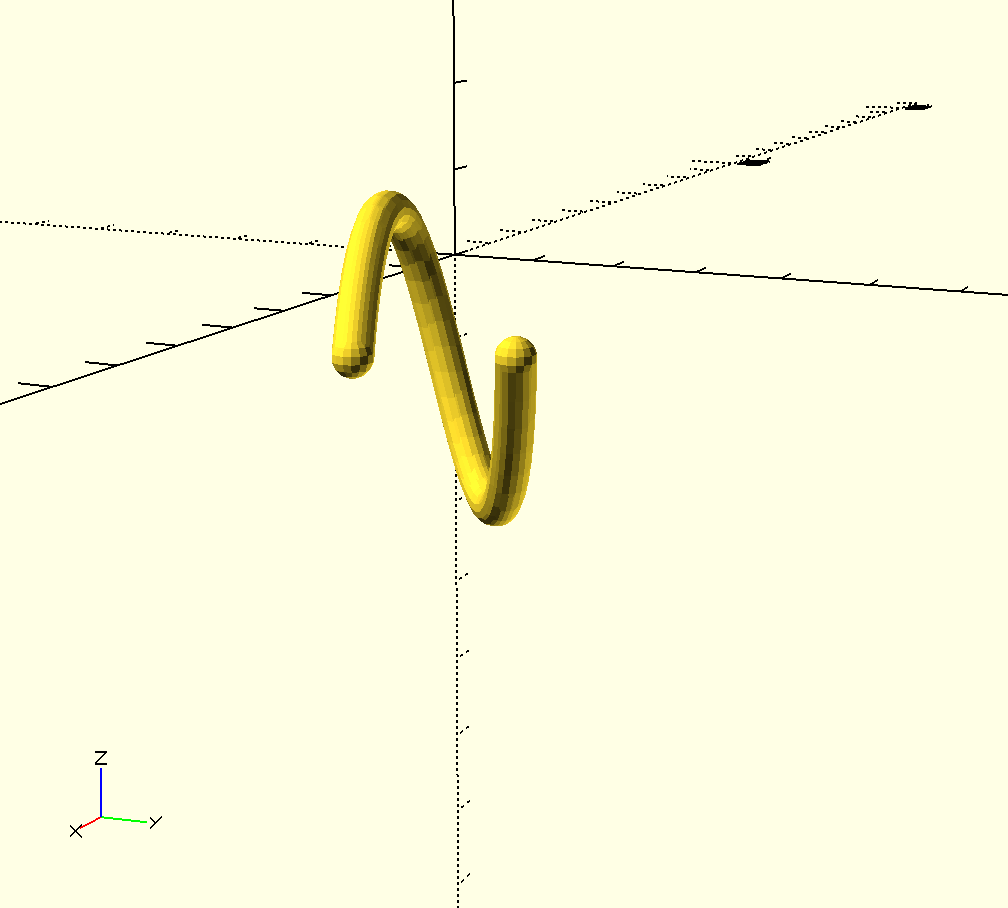} \\[0.02\textwidth]
  \includegraphics[width=0.49\textwidth]{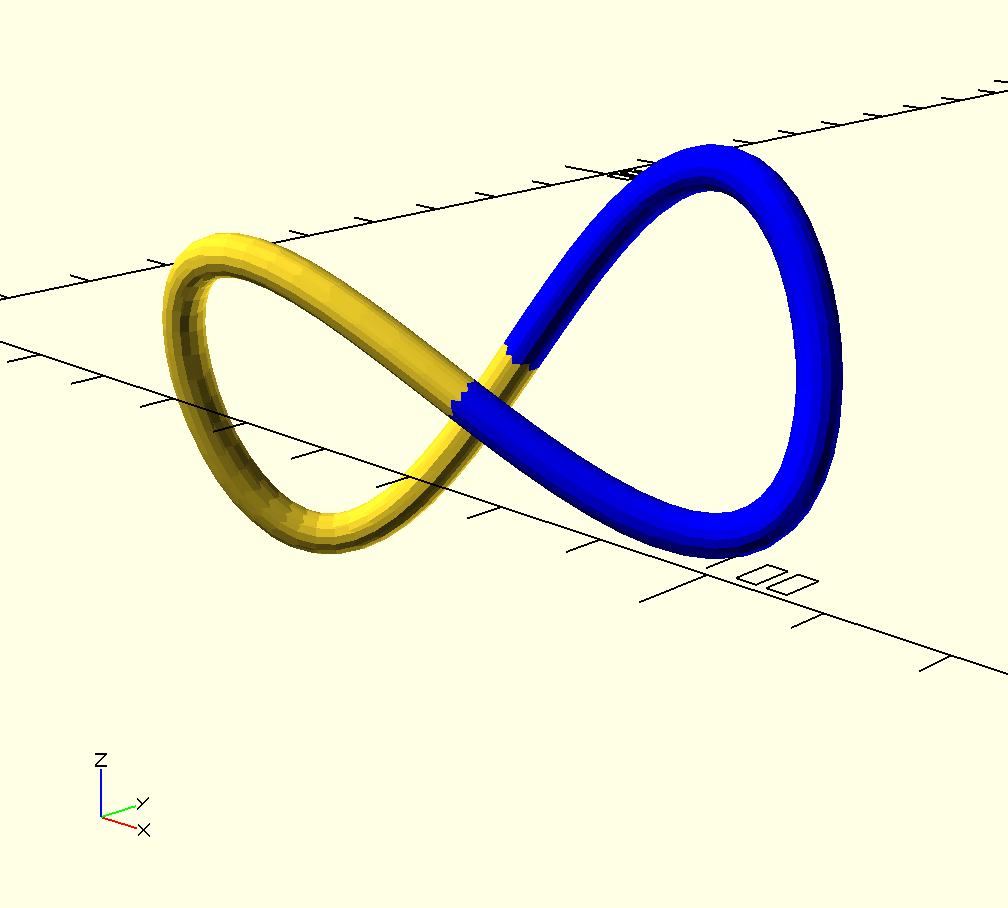}\hfill
  \includegraphics[width=0.49\textwidth]{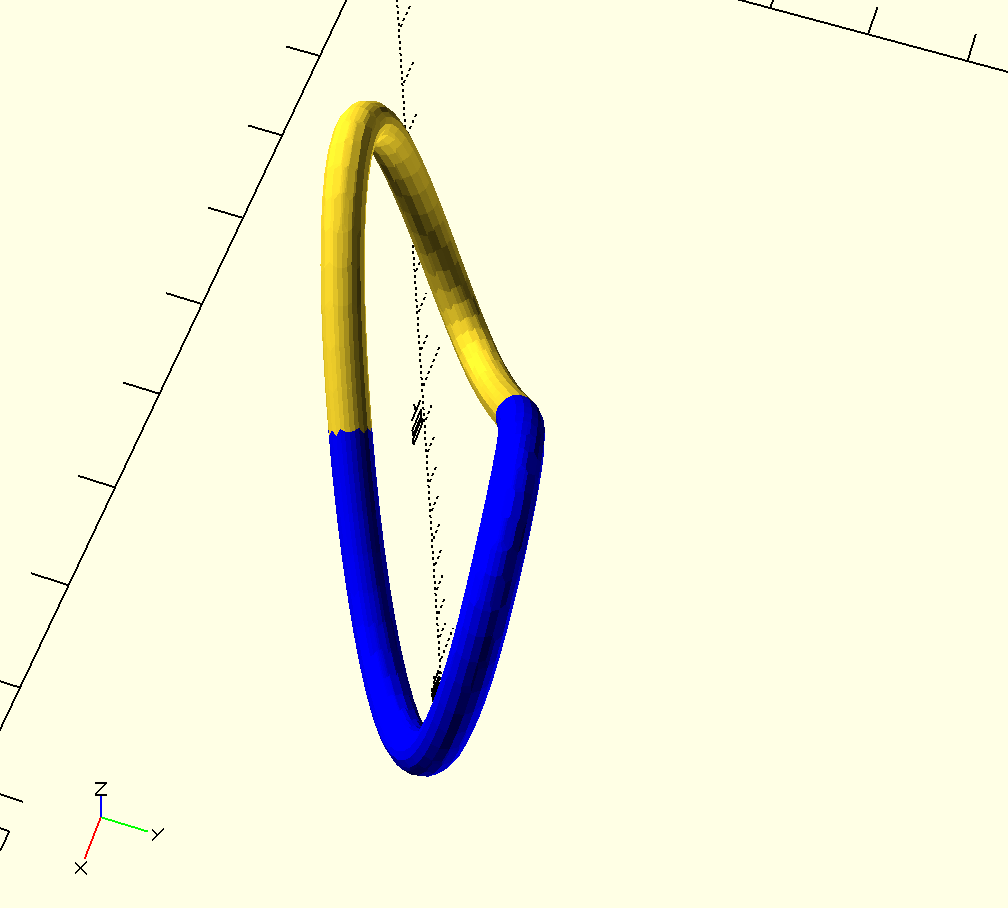}
\end{center}
  \caption{Top row: 
Half-lemniscate with oscillation to use as a building block for the transversal fibers in the Klein bottle model. 
The curling ensures that the transversal ribs of the wireframe do not selfintersect. 
\newline
Bottom row:
As copies of this half-lemniscate are placed along a Möbius strip path, the rotation of the opposite copy produces a seamless connection between the pieces.
The two copies that line up connect with each other to form a non-selfintersecting curve with a lemniscate projection.}
  \label{fig:half-lemniscate}
\end{figure}

\textbf{Klein Bottle}
For the Klein bottle, I had to take more care to get the dodged self-intersections right.
First, I used a helper function that parametrizes a point along half a lemniscate, with some oscillation orthogonal to the curve to pull apart the end points.
The oscillation can be seen in Figure~\ref{fig:half-lemniscate}.
To get this oscillation, I started out with the Lemniscate of Gerono parametrization
\[
u\mapsto (\cos u, \sin u \cdot \cos u)
\]

Next, I added a third dimension to the parametrization, with an offset cosine curve
\[
u\mapsto \left(\cos u, \sin u \cdot \cos u, \alpha\cos\left(\frac{\pi}{2}+u\right)\right)
\]
The $\pi/2$ offset was found by experimentation.
A choice of $\alpha=0.25$ turned out to give an aesthetically pleasing wobble -- where adjacent transversal curves space evenly along the surface.

These dodged half lemniscate were then translated and rotated to place them along a twisting Möbius strip shaped path.
This way, one full period after each curve was placed, its companion is placed rotated by $\pi$ around two different axes, placing it for a perfect gluing, as seen in Figure~\ref{fig:half-lemniscate}.
The resulting OpenSCAD code can be seein in Listing~\ref{lst:klein}.

\begin{lstfloat}
\begin{lstlisting}
function half_lemniscate(alpha, ampl) = 
  let(beta = 90+180*alpha)
  [cos(beta), sin(beta)*cos(beta), ampl*cos(90+beta)];

function klein(i,j) = 
  let(a_i = 360*i/lat_ribs)
  let(a_j = 360*j/long_ribs)
  let(pt = outer_radius*[1,0,0] + 
    inner_radius*rotate_matrix(0,0,a_i/2)*
    half_lemniscate(a_j/360, 0.25))
  let(pt = rotate_matrix(90,0,0)*pt)
  let(pt = rotate_matrix(0,0,a_i)*pt)
  pt;
\end{lstlisting}
\caption{Parametrization code for the Klein bottle immersion.
The \lstinline|ampl*cos(90+beta)| component makes the curve oscillate slightly, so that the self-intersection of the Klein bottle is highlighted by the two intersecting grid lines interweaving.}
\label{lst:klein}
\end{lstfloat}

\textbf{Projective Plane -- Roman surface}
For the projective plane, I chose the Roman surface parametrization.
This surface, also known as the Steiner surface is a mapping of the projective plane to $\mathbb{R}^3$ with tetrahedral symmetry with six singular points.

The Roman surface can be constructed by taking the unit sphere in $\mathbb{R}^3$ and mapping it using
\[
(x,y,z)\mapsto(yz,xz,xy)
\]

In my OpenSCAD code in Listing~\ref{lst:roman}, this mapping is paired with a parametrization of the sphere as
\[
(u,v)\mapsto( R\cos u\cdot\cos v, R\cos u\cdot\sin v, R\sin u)
\]

\begin{lstfloat}
  \begin{lstlisting}
function roman(i,j) = 
  let(u = i*360/lat_ribs)
  let(v = j*180/long_ribs)
  let(cu = cos(u)) 
  let(cv = cos(v))
  let(sv = sin(v))
  let(x=outer_radius*cu*cv)
  let(y=outer_radius*cu*sv)
  let(z=outer_radius*su)
  [y*z, x*z, x*y];
  \end{lstlisting}
\caption{Parametrization code for the Projective plane immersion.
For the Roman surface, we use the Veronese embedding directly, first constructing $x$, $y$ and $z$ on a sphere, and then emitting $(yz, xz, xy)$ for the parametrization itself.}
\label{lst:roman}
\end{lstfloat}

\textbf{Wireframe Construction}
To construct the wireframe representation, the same core loop is used in all three cases.
Writing \lstinline{pos} for the parametrizing function, and \lstinline{pipe} for the function to construct a capsule, the core loop can be seen in Listing \ref{lst:core}.
The loop goes through all points in a rectangular grid in the $u-v$ plane, and for each point uses it is as the origin for two segments -- one from $(u,v)$ to $(u+1,v)$ and then one from $(u,v)$ to $(u,v+1)$.
Along each of these two segments, capsules are placed along short steps using \lstinline{pos} to calculate coordinates for their placements.

\begin{lstfloat}
  \begin{lstlisting}
for(i=[0:1:2*lat_ribs]) {
  for(j=[0:1:long_ribs]) {
    for(x=[i:1/outer_density:i+1]) {
      pipe(pos(x,j), pos(x+1/outer_density,j), thickness);
    }
    for(y=[j:1/inner_density:j+1]) {
      pipe(pos(i,y), pos(i,y+1/inner_density), thickness);
    }
  }
}    
  \end{lstlisting}
\caption{Core loop constructing the wireframe.
Here, we use \lstinline|pipe| for the function that constructs a capsule, and \lstinline|pos| for the parametrizing function.
The code traverses the parameter grid, and for each grid point has two loops -- one to construct the geodesic edge $(i,j) \to (i,j+1)$ and one to construct the geodesic edge $(i,j) \to (j+1,j)$.
}
\label{lst:core}
\end{lstfloat}

\printbibliography
\end{document}